\documentclass[12pt]{amsart}

\usepackage[]{amsmath, amsthm, amsfonts}
\usepackage{amscd}
\usepackage{amssymb}            
\usepackage[mathcal]{euscript}

\usepackage[matrix,arrow,curve,frame]{xy}       % XY-pic diagram pac
\xymatrixcolsep{1.9pc}                          % Adjust size of diagrams.
\xymatrixrowsep{1.9pc}
\newdir{ >}{{}*!/-5pt/\dir{>}}                  % Make better tailed arrows

\newtheorem{theorem}{Theorem}[section]
\newtheorem{proposition}[theorem]{Proposition}
\newtheorem{lemma}[theorem]{Lemma}

\theoremstyle{definition}
\newtheorem{definition}[theorem]{Definition}

\newtheorem{notation}[theorem]{Notation}

\theoremstyle{remark}
\newtheorem{remark}[theorem]{Remark}

\numberwithin{equation}{section}

\DeclareMathOperator{\hocolim}{hocolim}
\DeclareMathOperator{\colim}{colim}
\DeclareMathOperator{\Ne}{N}

\newcommand{\Coprod}{\,{\scriptstyle\coprod}\,}
\newcommand{\tcoprod}{\,{\textstyle\coprod}\,}
\newcommand{\Sum}{\,{\scriptstyle\sum}\,}

\newcommand{\A}{{\mathcal A}}
\newcommand{\B}{{\mathcal B}}
\newcommand{\C}{{\mathcal C}}

\newcommand{\I}{{\mathcal I}}
\newcommand{\N}{{\mathcal N}}

\renewcommand{\S}{{\mathcal S}}
\newcommand{\x}{{\mathbf x}}
\renewcommand{\O}{{\mathcal O}}
\newcommand{\la}{\langle}
\newcommand{\ra}{\rangle}
\newcommand{\ov}{\overline}

\newlength{\labwidth}
\newcommand{\labarrow}[1]{
\settowidth{\labwidth}{$\scriptstyle \;\; #1 \;\;$}
\stackrel{#1}{\smash{\hbox to \labwidth{\rightarrowfill}} 
\vphantom{\longrightarrow}}
}

\begin{document}

\title{Multivariable cochain operations and little $n$-cubes.}
\author[J.E. McClure]{James E. McClure} 
\address{
Department of Mathematics, Purdue University  \\
150 N. University Street \\
West Lafayette, IN  47907-2067}
\email{mcclure@math.purdue.edu}
\thanks{The first author was partially supported by NSF grant DMS-9971953. He
thanks the Lord for making his work possible.}

\author[J.H. Smith]{Jeffrey H. Smith}
\address{
Department of Mathematics, Purdue University  \\
150 N. University Street \\
West Lafayette, IN  47907-2067}
\email{jhs@math.purdue.edu}
\thanks{The second author was partially supported by NSF grant DMS-9971953.}


\date{June 7, 2001 and, in revised form, June 27, 2002.}

\begin{abstract}
In this paper we construct a small $E_\infty$ chain operad $\S$ which acts
naturally on the normalized cochains $S^*X$ of a topological space.  We also
construct, for each $n$,  a suboperad $\S_n$ which is quasi-isomorphic to the
normalized singular chains of the little $n$-cubes operad.  The case $n=2$
leads to a substantial simplification of our earlier proof of Deligne's
Hochschild cohomology conjecture.
\end{abstract}

\maketitle

\section{Introduction.}

This paper has two goals.  The first (see Theorem \ref{revS2} and
Remark \ref{revM5}(a)) is to 
construct a small $E_\infty$ chain operad $\S$ which acts naturally on the 
normalized cochains $S^*X$ of a topological space $X$.  This is of interest 
in view of a theorem of Mandell \cite[page 44]{Mandell} which states that if 
$\O$ is any $E_\infty$ chain operad over $\ov{\mathbb F}_p$ (the algebraic 
closure of the field with $p$ elements) which acts naturally on $S^*X\otimes 
\ov{\mathbb F}_p$ then the homotopy category of connected $p$-complete 
nilpotent spaces of finite type imbeds in the homotopy category of $\O$ 
algebras; tensoring our operad with $\ov{\mathbf F}_p$ gives an operad to
which Mandell's theorem applies.  Operads which act naturally on the 
normalized cochains of spaces were already known to exist \cite{Dold1}, 
\cite{HinichSchechtman}, but the known examples are uncountably generated 
and the structure maps are hard to describe explicitly,
whereas our operad is of finite type and has a perspicuous description (see
Definition \ref{persp} and Propositions \ref{diff}, \ref{perm} and \ref{gamma}).

The second goal (see Theorem \ref{big}) is to give, for each $n$, a 
suboperad $\S_n$ of $\S$ which
is quasi-isomorphic, in the category of chain operads over the integers, to the
normalized singular chain operad of the little $n$-cubes operad (this is the
smallest and simplest known chain model for the little $n$-cubes).  
One reason this is interesting is that the operad
$\S_2$ acts naturally on the normalized Hochschild cohomology complex of an
associative ring (see Theorem \ref{De}), so we obtain a solution of Deligne's 
Hochschild cohomology conjecture \cite{Deligne}.   Other solutions of 
Deligne's Hochschild cohomology conjecture are known 
(\cite{K1,K2,MS,T1,T2,V}) 
but the one we give here is significantly shorter than the previously known 
solutions.
Also, our solution (like those in \cite{MS,K2}) is valid over the integers and
not only in characteristic 0.
We expect that Theorem \ref{big} will lead to one way of proving 
the generalized version of Deligne's conjecture proposed by Kontsevich in
\cite{K1}.

The operad $\S_2$ is isomorphic to the operad $\C$ considered in our previous
paper \cite{MS}, but the description we give here is more explicit and the
proof that $\S_2$ is quasi-isomorphic to the normalized singular chains of 
the little 2-cubes is much simpler than the corresponding proof in \cite{MS}.

The organization is as follows: in Section \ref{sec2} we define the chain 
operad $\S$ and state Theorem \ref{revS2} and Propositions \ref{diff},
\ref{perm} and \ref{gamma}.
In Section \ref{sec3} we define the chain operads $\S_n$
and state Theorem \ref{big}.  
In Section \ref{sec3a} we describe the action of $\S_2$ on the Hochschild
complex.  Section \ref{sec7} gives the proof of Theorem \ref{big}, using a
method due to Berger \cite{Berger}, and Sections \ref{sec4}, \ref{sec5} and 
\ref{sec6} give the proofs of Propositions \ref{diff}, \ref{perm} and
\ref{gamma}.

The results of section \ref{sec2} have been discovered independently 
(but with a different choice of signs) in joint work of Berger and Fresse 
\cite{BergerFresse}.  After learning Theorem \ref{big} from us, they have 
found a different proof of it using their methods.

\section{The chain operad $\S$.}
\label{sec2}

Our goal in this section is to define a family of cochain operations which
generate a small $E_\infty$ chain operad $\S$ and to give the formulas which
describe the structure of $\S$.

First we need to specify what we mean by cochains.
We write $S_* X$ for the normalized singular chains of a space $X$ (that is, 
the singular chains modulo degeneracies; the reason for using normalized 
chains is given in the proof of Lemma \ref{revS1}).  

\begin{definition}
Let ${\mathbb Z}[0]$ denote ${\mathbb Z}$
considered as a chain complex concentrated in dimension 0.  The {\it normalized
singular cochain complex} $S^* X$ is the cochain complex Hom$(S_* X,{\mathbb 
Z}[0])$.
\end{definition}

Note that with this definition the coboundary $d$ satisfies
\[
d(x)=-(-1)^{|x|}x\circ \partial
\]
which differs by a sign from the usual coboundary (see \cite[Remark
VI.10.28]{Dold}).

Before giving the formal definitions of the cochain operations we are
interested in, we begin with an informal discussion;
in particular we denote all signs by $\pm$ for the time being.

\begin{definition}
If $\sigma:\Delta^r \to X$, and if 
$a_0,\ldots,a_p\in \{0,\ldots,r\}$, we write
\[
\sigma(a_0,\ldots,a_p):\Delta^p\to\Delta^r
\]
for the affine map which takes vertex $i$ of $\Delta^p$ to vertex $a_i$ of
$\Delta^r$.
\end{definition}

Our first example of a cochain operation is the chain-level cup product.  
If 
$x\in S^p X$, $y\in S^q X$, then $x\smallsmile y\in S^{p+q} X$ is defined 
by 
\[
(x\smallsmile y)(\sigma)=x(\sigma(0,\ldots,p))\cdot y(\sigma(p,\ldots,p+q))
\]
where $\sigma:\Delta^{p+q}\to X$.

Next we recall Steenrod's original definition of the cup-$1$ product 
\cite[bottom of page 293]{Steenrod}.
If $x\in S^p X$, $y\in S^q X$, then
$x\smallsmile_1 y\in S^{p+q-1} X$ is defined by
\[
(x\smallsmile_1 y)(\sigma)=\sum_{i<j} \pm x(\sigma(0,\ldots,i,j,\ldots,p+q-1))\cdot
y(\sigma(i,\ldots,j))
\]
where $\sigma:\Delta^{p+q-1}\to X$.  (Some of the terms in this formula
involve evaluating a cochain on a simplex whose degree does not match that of
the cochain; all such terms are defined to be 0.)

To describe the higher cup-$i$ products we need a definition. 

\begin{definition}
Let $T$ be a finite totally ordered set.  An {\it overlapping partition} $\A$ 
of $T$ with
$m$ pieces is a collection of subsets $A_1,\ldots,A_m$ of $T$ with the 
following properties:

(a) If $j<j'$ then each element of $A_j$ is $\leq$ each element of $A_{j'}$.

(b) $A_j\cap A_{j+1}$ has exactly one element for each $j<m$.
\end{definition}

Thus consecutive pieces of the overlapping partition are required to have
exactly one element in common.
For example, the collection $\{\{0,1,2\},\{2,3\},\{3\},\{3,4,5\}\}$ is an
overlapping partition of $\{0,1,2,3,4,5\}$.

\begin{remark}
\label{revM2}
For later use we record the elementary fact that an overlapping partition 
$A_1,\ldots,A_m$ is completely determined by the $m-1$ elements $A_1\cap A_2,
\ldots,A_{m-1}\cap A_m$, which we will refer to as the {\it overlap points}.
\end{remark}

Now Steenrod's original definition of the cup-$i$ product \cite[top of page
294]{Steenrod} is
\[
(x\smallsmile_i y)(\sigma)=\sum \pm 
x(\sigma(A_1 \Coprod A_3\Coprod\ldots))\cdot
y(\sigma(A_2\Coprod A_4\Coprod \ldots))
\]
where the sum is over all overlapping partitions of
$\{0,\ldots,p+q-i\}$ with $i+2$ pieces.  Here $A_1\Coprod A_3\Coprod\ldots$ 
and $A_2\Coprod A_4\Coprod \ldots$ denote the {\it
disjoint} unions, that is, repetitions are not eliminated; note that 
all terms for which $A_1\Coprod A_3\Coprod\ldots$ or 
$A_2\Coprod A_4\Coprod \ldots$ have repetitions will be zero, because the
corresponding simplices will be degenerate.

There is a useful piece of notation (due to Benson \cite[page 147]{Benson}  
and Milgram,
and rediscovered by us) which clarifies the situation.  We write $12(x,y)$ for
$x\smallsmile y$, $121(x,y)$ for 
$x\smallsmile_1 y$, $1212(x,y)$ for 
$x\smallsmile_2 y$, and so on.  We also write $21(x,y)$ for $y\smallsmile x$,
$212(x,y)$ for $y\smallsmile_1 x$ and so on. The idea is that to evaluate (for
example) $12121(x,y)$ on a simplex $\sigma:\Delta^r\to X$ we sum over all
overlapping partitions of $\{0,\ldots,r\}$ with 5 pieces (because the sequence 
12121 has 5 entries) and the term corresponding to a partition 
$A_1,A_2,A_3,A_4,A_5$ is 
\[
x(\sigma(A_1\Coprod A_3\Coprod A_5))\cdot
y(\sigma(A_2\Coprod A_4));
\]
here $A_1,A_3,A_5$ go with $x$ because the 1st, 3rd and 5th entries of the
sequence 12121 are 1's and $x$ is the 1st variable in the symbol $12121(x,y)$.

It will now be obvious that one can use a similar idea to define multilinear 
cochain operations with more than 2 variables (this fact was first noticed by 
Benson \cite[page 147]{Benson}, who used it to give a definition of the 
odd-primary Steenrod operations in the cohomology of groups, and by 
Milgram, and later independently by us).  For 
example, it is natural to interpret the symbol $12312(x,y,z)$ by
%\begin{multline}
%\notag
\[
12312(x,y,z)(\sigma)=\sum \pm x(\sigma(A_1\Coprod A_4))\cdot 
y(\sigma(A_2\Coprod A_5))
\cdot z(\sigma(A_3))
\]
%\end{multline}

Similarly one can define a multilinear cochain operation in $k$ variables for
any finite sequence with entries $1,\ldots,k$ in which each of the numbers
$1,\ldots,k$ is used at least once.
We will call multilinear cochain operations of this type {\it sequence
operations}.

Now we turn to the formal definitions.

\begin{notation}
For $k\geq 1$ let $\bar{k}$ be the set $\{1,\ldots,k\}$, and let
$\bar{0}=\emptyset$.
\end{notation}

It will be convenient to think of sequences of length $m$ with entries in the 
set $\bar{k}$ as functions $\bar{m}\to\bar{k}$.  The sequences in which each 
of the numbers $1,\ldots,k$ is used at least once correspond to 
surjections $\bar{m}\to\bar{k}$.  

Our next three definitions are needed in order to specify the signs in our
formulas.

\begin{notation}
The cardinality of a set $A$ will be denoted $|A|$, and $||A||$ will denote
$|A|-1$.
\end{notation}

\begin{definition}
\label{tau}
For each surjection $f:\bar{m}\to\bar{k}$ define $\tau_f:\bar{m}\to\bar{m}$ by
\[
\tau_f(j)=|\{j'\in\bar{m}\,:\, f(j')<f(j) \mbox{ or $f(j')=f(j)$ and $j'\leq 
j$}\}|
\]
\end{definition}

\begin{remark}
Any surjection $f:\bar{m}\to\bar{k}$ can be factored (usually in more than one
way) as a permutation of $\bar{m}$ followed by an order-preserving map
$\bar{m}\to\bar{k}$.  Among permutations of $\bar{m}$ that arise in this 
way, $\tau$ is the only one that is order preserving on each set $f^{-1}(i)$.
\end{remark}

\begin{definition}
\label{revT3}
Let $T$ be a finite totally ordered set.
For each surjection $f:\bar{m}\to\bar{k}$ and each overlapping partition $\A$ 
of $T$ with $m$ pieces, define
\[
\epsilon(f,\A)= 
\underset{f(j)>f(j')}{\underset{j<j'}\sum}
||A_j||\,||A_{j'}||
+
\sum_{j=1}^m ||A_j||(\tau_f(j)-f(j))
\]
\end{definition}

Now we can define the sequence operation associated to a surjection
$f:\bar{m}\to\bar{k}$.  We write $S_*X$ for the normalized singular chains of
$X$.

\begin{definition}
\label{la}
(a)
Given a surjection $f:\bar{m}\to\bar{k}$ and 
$\sigma:\Delta^p\to X$ 
define $\sigma[f]\in (S_*X)^{\otimes k}$ by
\[
\sigma[f]
=
\sum_\A (-1)^{\epsilon(f,\A)} 
\bigotimes_{i=1}^k \,
\sigma\Bigl(\,{\mathop{{\textstyle{\coprod}}}_{f(j)=i}} \, A_j\Bigr)
\]
Here $\A$ runs through the overlapping partitions of $\{0,\ldots,p\}$ with 
$m$ pieces. 

(b)
Given a surjection $f:\bar{m}\to\bar{k}$, define a natural transformation
\[
\la f \ra: (S^*X)^{\otimes k} \to S^* X
\]
by
\[
\la f \ra(x_1\otimes \ldots\otimes 
x_k)(\sigma)=(-1)^{(m-k)}(x_1\otimes\cdots\otimes x_k)
(\sigma[f]) 
\]
\end{definition}

\begin{remark}
\label{Fri1}
(a) When $f:\bar{m}\to \bar{2}$ is the function corresponding to the sequence
$1212\cdots$ with $m$ entries then $\la f\ra$ is (up to sign) the same as
Steenrod's $\smallsmile_{m-2}$.  

(b) In particular, when $f$ corresponds to 12 the
operation $\la f\ra (x_1\otimes x_2)$ is the same as Dold's definition 
\cite[page 222]{Dold}
of the cup product $x_1\smallsmile x_2$,
which differs from the usual definition by the sign $(-1)^{|x_1||x_2|}$.

(c) For all $f$, $\la f\ra$ lowers the total degree by $m-k$.

(d) To explain the sign in part (b) of Definition \ref{la} we first note that 
the degree of $\sigma$ must be $(|x_1|+\ldots+|x_k|)-(m-k)$ in order for the
left-hand side to be nonzero.  Now
$\la f \ra$, which has degree $m-k$, is being moved past $x_1,\ldots,x_k$ and
$\sigma$, so the sign is $(-1)^{(m-k)^2}=(-1)^{m-k}$.

(e) There are other ways to choose the sign in part (a) of Definition
\ref{la},
but the choice we have given seems to lead to the simplest signs 
in other formulas (such as those in Propositions \ref{diff}, \ref{perm} and
\ref{gamma}) that depend on this one.  Our choice of signs is motivated by the
cosimplicial point of view developed in \cite{MS3}.

(f) Here is a conceptual way to keep track of the signs: we associate to $f$ 
the product of simplices $\prod \Delta^{||f^{-1}(i)||}$, and we represent the
$||f^{-1}(i)||$-dimensional simplex $\Delta^{||f^{-1}(i)||}$ by a sequence of
$||f^{-1}(i)||$ one-dimensional objects
(each denoted by $*$). The sign in part (a) of Definition \ref{la} comes from 
permuting the
simplices $\sigma(A_j)$ and then moving the $*$'s corresponding to 
$\Delta^{||f^{-1}(i)||}$, in order, into the positions
between the pieces of $\sigma(\coprod_{f(j)=i} \, A_j)$.
\end{remark}

Next we observe that $\la f \ra$ 
will be identically zero for certain $f$.

\begin{lemma} 
\label{revS1}
If $f(l)=f(l+1)$ for some $l\in \bar{m}$ then $\la f\ra=0$
\end{lemma}

\begin{proof}
It suffices to show that $\sigma[f]=0$ for all $\sigma: \Delta^p\to X$.  
Let $i_0=f(l)=f(l+1)$.  If 
$\A$ is any
overlapping partition of $\{0,\ldots,p\}$ with $m$ pieces, then (since $A_l$ 
and $A_{l+1}$ have an element in common) 
${\mathop{{\textstyle{\coprod}}}_{f(j)=i_0}} \, A_j$ will have a
repeated entry and therefore the simplex
$\sigma({\mathop{{\textstyle{\coprod}}}_{f(j)=i_0}} \, A_j)$ will be
degenerate, which means it will represent 0 in the normalized chains $S_*X$.
\end{proof}

This motivates:

\begin{definition}
\label{ndg}
A function $f:\bar{m}\rightarrow\bar{k}$ is {\it nondegenerate} if it is
surjective and $f(j)\neq f(j+1)$ for $1\leq j< m$.  Otherwise $f$ is {\it
degenerate}.
\end{definition}

In particular, the identity map of $\bar{0}$ is nondegenerate.  The sequence
corresponding to a nondegenerate function has adjacent entries distinct.  

\begin{definition}
\label{persp}
Let 
$\S(k)$ be the graded abelian group freely generated by the 
maps $f:\bar{m}\to\bar{k}$ (where $f$ is assigned the degree 
$m-k$) modulo the subgroup generated by the degenerate maps. 
\end{definition}

The symbol $\S$ stands for ``sequence operad.''
Note that $\S(k)$ is freely generated by the nondegenerate 
$f:\bar{m}\to\bar{k}$.

Now let $\N(k)$ be the graded abelian group of natural transformations 
\[
(S^* X)^{\otimes k}\to S^* X
\]
where a transformation $\nu$ is assigned degree $n$ if it lowers total degree 
by $n$.  Definition \ref{la} gives a homomomorphism $\S(k)\to \N(k)$ which is 
easily seen to be a monomorphism.  Each $\N(k)$ is a chain complex, with 
differential
\[
\partial \nu =d \circ \nu -(-1)^{|\nu|} \nu\circ d
\]
where $d$ is the coboundary of $S^*X$ (resp.\ $(S^* X)^{\otimes k}$).
Moreover, the collection $\N$ of all the $\N(k)$ with $k\geq 0$ is a chain
operad, where the action of the symmetric group on $\N(k)$ and the 
multivariable composition operations are the obvious ones. ($\N$ should be 
thought of as the endomorphism operad of the functor $S^*$.)

Let $\S$ be the collection of all $\S(k)$ with $k\geq 0$.
Our main theorem in this section is:

\begin{theorem}
\label{revS2}
(a) $\S(k)$ is a subcomplex of $\N(k)$.

(b) $\S$ is a sub-chain-operad of $\N$.

(c) Each $\S(k)$ has the homology of a point.
\end{theorem} 

\begin{remark}
\label{revM5}
(a)
Part (c) of the theorem says that $\S$ is an $E_\infty$ chain operad.  Thus
$\S$ is an $E_\infty$ chain operad which (by definition) acts naturally on 
$S^*X$.  Propositions \ref{diff}, \ref{perm} and \ref{gamma} below give
explicit formulas for the structure of $\S$.

(b)
A right coalgebra over an operad 
$\O$ is a chain complex $C$
together with maps
\[
\delta_k:C\otimes\O(k)\to C^{\otimes k}
\]
which are consistent with the symmetric group action, composition and unit of
$\O$.
The maps $S_*X\otimes \S(k)\to (S_*X)^{\otimes k}$ which take
$\sigma\otimes f$ to $\sigma[f]$ 
give $S_*X$ a natural structure of right coalgebra over $\S$.

(c) Part (a) of the theorem is stated on page 147 of \cite{Benson} without
a proof.  The same page states incorrectly that $\S(k)$ is not contractible,
but it gives the correct chain homotopy needed to prove (c).

(d) Ezra Getzler has informed us that during the period from 1995 to 1997, 
after learning of the sequence operations from Milgram, he proved part (c) of 
Theorem \ref{revS2} and noticed parts (a) and (b) without proving them.
\end{remark}

The proof of Theorem \ref{revS2} will occupy the rest of this section, with
some calculations deferred to later sections.

To prove part (a) we need to show that the differential of a sequence operation
is a ${\mathbb Z}$-linear combination of sequence operations.  
To state the
precise result we need a definition.

\begin{definition}
\label{tues1}
Given $f:\bar{m}\to\bar{k}$ and a subset $S$ of $\bar{m}$, let $s=|S|$ and 
write $f_S$ for the composite
\[
\bar{s}\to S \hookrightarrow \bar{m}\labarrow{f}\bar{k},
\]
where the first map is the unique order-preserving bijection.
\end{definition}

\begin{proposition}
\label{diff}
Let $f:\bar{m}\to\bar{k}$ be nondegenerate. Then
\[
\partial\la f\ra=
\sum_{j=0}^m \, (-1)^{\tau_f(j)-f(j)}\, \la f_{\bar{m}-\{j\}}\ra
\]
\end{proposition}

(Note that, by Lemma \ref{revS1}, the terms for which 
$f_{\bar{m}-\{j\}}$ is degenerate are zero.)

For example, if we write $f$'s as sequences we have
\[
\partial\la 12312\ra=\la 2312\ra -\la 1232\ra -\la 1312\ra +\la 1231\ra 
\]
and
\[
\partial\la 123121\ra =\la 23121\ra -\la 12321\ra +\la 12312\ra +\la 13121\ra 
\]
The pattern is that we delete each 1, then each 2, etc., and the signs
alternate, except that the last term obtained by deleting $i$ has the same sign
as the first term obtained by deleting $i+1$; if a term has an adjacent pair 
equal or if it doesn't contain all the numbers from 1 to $k$ it is 
zero.  (A conceptual way to remember the signs is that they are the same as 
the signs in the cellular chain complex of
$\prod\Delta^{||f^{-1}(i)||}$; cf.\ Remark \ref{Fri1}(f).)

Here is an outline of the proof of Proposition \ref{diff}; the complete proof 
is given in Section \ref{sec5}.  
By definition $(\partial\la 
f\ra)(x_1\otimes \ldots\otimes x_k)(\sigma)$ is
\[
\pm \la f\ra(x_1\otimes \ldots\otimes x_k)(\partial \sigma) 
\pm \la f \ra(d(x_1\otimes \ldots\otimes x_k))(\sigma)
\]
and according to Definition \ref{la}(b) this can be rewritten
\begin{equation}
\label{reve1}
\pm (x_1\otimes \ldots\otimes x_k)((\partial \sigma)[f])
\pm (x_1\otimes \ldots\otimes x_k)(\partial(\sigma[f]))
\end{equation}
so the question is what the difference is between 
$(\partial \sigma)[f]$ and $\partial(\sigma[f])$; the statement of Proposition
\ref{diff} says that this difference should be
a sum of terms of the form $\pm \sigma[f_{\bar{m}-\{j\}}]$.
When we evaluate $\partial(\sigma[f])$ using Definition \ref{la}(a), we get
a sum indexed over overlapping partitions $\A$ and, for each such
partition, a sum in which the elements of the sets 
$\coprod_{f(j)=i} \, A_j$ are deleted in turn.  Most of the terms in this
double sum also occur in 
$(\partial \sigma)[f]$ and so they cancel in (\ref{reve1}) (assuming
that the signs work correctly, which they do).  If one of the sets $A_j$ has a
single element, then the term in which it is deleted does not occur in 
$(\partial \sigma)[f]$;
instead it is a term in $\sigma[f_{\bar{m}-\{j\}}]$.  The only other
possibility is that $|A_j|\geq 2$ and the element in $A_j$ which is deleted 
is at the beginning or end of $A_j$; terms of this form occur twice each
in $\partial(\sigma[f])$, and
when the signs are taken into account they cancel.

We now turn to part (b) of Theorem \ref{revS2}.  We must show three things:
that the identity element of $\N(1)$ is in $\S(1)$, that $\S(k)$ is preserved
by the permutation action on $\N(k)$, and that $\S$ is closed under
multivariable composition.  The first of these is easy: the identity element in
$\N(1)$ is the sequence operation induced by the identity map
$\bar{1}\to\bar{1}$.

Next we must calculate the effect of permutations on sequence operations.
Inspection of Definition \ref{la} shows at once that the composite of a
permutation and a sequence operation is (up to sign) another sequence
operation.  Here is the precise formula:

\begin{proposition}
\label{perm}
Let $\rho\in \Sigma_k$ and let $f:\bar{m}\to\bar{k}$ be nondegenerate.
Then 
\[
\la f\ra\rho=(-1)^{\zeta(f,\rho)}\la \rho^{-1}\circ f\ra
\]
where the sign is given by
\[
\zeta(f,\rho) =\sum
||f^{-1}(i)||\,||f^{-1}(i')||
\]
with the sum taken over all pairs for which $i<i'$ and
$\rho^{-1}(i)>\rho^{-1}(i')$.
\end{proposition}

A conceptual way to remember the sign is that it is the same as the sign coming
from the permutation of the factors in $\prod
\Delta^{||f^{-1}(i)||}$; cf.\ Remark \ref{Fri1}(f).

We defer the proof of Proposition \ref{perm} to Section \ref{sec4}.

To complete the proof of part (b) of Theorem \ref{revS2} we need to
show that the multivariable composite of sequence operations is a 
${\mathbb Z}$-linear combination of sequence operations.  The formula may be
found in Proposition \ref{gamma}, but it is complicated so
we begin with some motivation.

Let $f:\bar{m}\to\bar{k}$ and let $g_i:\overline{m_i}\to\overline{j_i}$ for
$1\leq i\leq k$.  We want to describe the multivariable composition
$\la f\ra(\la g_1\ra, \ldots, \la g_k\ra)$.
First observe that
\begin{multline}
\notag
\la f\ra(\la g_1\ra,\ldots, \la g_k \ra)(x_1\otimes\ldots\otimes
x_{j_1+\cdots+j_k})\\
=
\pm \la f \ra(\la g_1\ra(x_1\otimes \cdots\otimes x_{j_1})\otimes 
\cdots\otimes \la g_k\ra(x_{j_1+\cdots+j_{k-1}+1}\otimes\cdots\otimes
x_{j_1+\cdots+j_k}))
\end{multline}
so we want to evaluate
\begin{equation}
\label{reve2}
\la f \ra(\la g_1\ra(x_1\otimes \cdots\otimes x_{j_1})\otimes
\cdots\otimes \la g_k\ra(x_{j_1+\cdots+j_{k-1}+1}\otimes\cdots\otimes
x_{j_1+\cdots+j_k}))(\sigma)
\end{equation}
To evaluate (\ref{reve2}) we apply
Definition \ref{la} twice and this leads to
a double sum, with the outer sum indexed by overlapping partitions 
$\A$ and the inner sum, for each $\A$, indexed by overlapping partitions of 
$\coprod_{f(j)=i} \, A_j$ for $1\leq i\leq k$.  This suggests the following
definition.

\begin{definition}
\label{revM3}
A {\it
double overlapping partition} $\mathbb A$ of $\{0,\ldots,p\}$ of type
$(f,m_1,\linebreak[0] \ldots,\linebreak[0] m_k)$
is an overlapping partition $\A$ of $\{0,\ldots,p\}$ with
$m$ pieces
together with, 
for each
$1\leq i\leq k$,
an overlapping partition $\B^i$ of
$\coprod_{f(j)=i} A_j$ with $m_i$
pieces.
\end{definition}

Thus (\ref{reve2}) is a sum indexed by the double overlapping partitions of 
type $(f,m_1,\linebreak[0] \ldots,\linebreak[0] m_k)$.  Next let us observe 
that the overlapping partition
$\A$ has $m-1$ ``overlap points'' (see Remark \ref{revM2})
and the overlapping partition $\B^i$ has $m_i-1$ overlap
points for each $i$; collecting these together we get $m-1+(\sum m_i)-k$ points,
and using these as overlap points we get a partition of $\{0,\ldots,p\}$ with 
$m-k+\sum m_i$ pieces: we will call this overlapping partition 
$\C$ and denote
its pieces by $C_1,\ldots,C_{m-k+\sum m_i}$.  Now in the 
expansion of the 
expression (\ref{reve2}) each of the $C_l$ is associated to one of the 
cochains $x_1,\ldots,x_{j_1+\cdots+j_k}$; this gives a surjection
\[
h:\overline{m-k+\Sum m_i}\to \overline{\Sum j_i}
\]
When all the terms in the expansion of (\ref{reve2}) which give rise to the
same $h$ are collected together, it is not hard (if signs are ignored) to see 
that they add to
\[
\la h\ra(x_1\otimes\cdots\otimes x_{j_1+\cdots+j_k})(\sigma)
\]
This shows that the multivariable composition $\la f\ra(\la g_1\ra, \ldots, 
\la g_k\ra)$ with which we began is indeed a $\mathbb Z$-linear combination of 
the sequence operations $\la h\ra$.

In order to say exactly which $h$'s occur for a given choice of
$f,g_1,\ldots,g_k$, 
and in order to specify the signs, and in order to be able to give 
efficient proofs, we need to reformulate the discussion just given in a more 
abstract way.

First we give an abstract description of the concept of overlapping partition.

\begin{definition}
\label{special}
Let $S$, $T$ and $U$ be finite totally ordered sets.
A diagram of the form
\[
\xymatrix{
S
&
T
\ar[l]_-{\alpha}
\ar[r]^-{\beta}
&
U
}
\]
is {\it special} if 

(i) $\alpha$ and $\beta$ are order-preserving epimorphisms, and

(ii) if $q<q'$ is an adjacent pair in $T$ then 
$\alpha(q)<\alpha(q')$ or
$\beta(q)<\beta(q')$ but not both
(this implies that $|T|=|S|+|U|-1$).
\end{definition}

\begin{lemma}
\label{part}
{\rm
A special diagram 
\[
\xymatrix{
\{0,\ldots,p\}
&
\overline{p+m}
\ar[l]_-{\alpha}
\ar[r]^-{\beta}
&
\bar{m}
}
\]
determines, and is determined by, an overlapping partition $\A$
of $\{0,\ldots,p\}$ with $m$ pieces, with $A_j=\alpha(\beta^{-1}(j))$.
}
\end{lemma}

The (easy) proof is left to the reader.

Our next three definitions are needed to specify the $\la h\ra$'s that arise in
evaluating a multivariable composition and the signs that go with them.

Fix $f:\bar{m}\rightarrow\bar{k}$
and nonnegative integers $m_1,\ldots,m_k$.
Let
\[
\chi:\coprod\overline{m_i}\rightarrow\bar{k}
\]
be the map which takes $\overline{m_i}$ to $i$.

\begin{definition}
\label{type}
A {\it diagram of type $(f,m_1,\ldots,m_k)$} is a commutative diagram of 
the form
\[
\xymatrix{
\overline{m-k+\Sum m_i}
\ar[r]^-{b}
\ar[d]_{a}
&
\coprod\overline{m_i}
\ar[d]_{\chi}
\\
\bar{m}
\ar[r]^{f}
&
\bar{k}
}
\]
such that 
$a$ is order preserving and 
for each $i\in\bar{k}$ the diagram
\[
\xymatrix{
f^{-1}(i)
&
b^{-1}(\overline{m_i})
\ar[r]^-{b}
\ar[l]_-{a}
&
\overline{m_i}
}
\]
is special.
\end{definition}

For example, a double overlapping partition of type $(f,m_1,\ldots,m_k)$ 
gives a diagram of type $(f,m_1,\ldots,m_k)$ as follows (see Definition 
\ref{revM3} and the paragraph which follows it for the notation): we define 
$a(l)=j$ if $C_l$ is contained in the $j$-th piece of the overlapping 
partition $\A$ and we define $b(l)=n\in \overline{m_i}$ if $C_l$ is contained 
in the $n$-th piece of the overlapping partition $\B^i$. 

\begin{definition}
Given a diagram of type $(f,m_1,\ldots,m_k)$ and $i\in\bar{k}$, let $\A^i$
be the overlapping partition of $f^{-1}(i)$ with $m_i$ pieces given by 
$\A^i_r=a(b^{-1}(r))$ for $r\in\ov{m_i}$.
\end{definition}

The diagram is in fact determined by the overlapping partitions $\A^i$; it 
expresses the
information contained in these overlapping partitions in a convenient form.

\begin{definition}
\label{HD}
Let $f:\bar{m}\to\bar{k}$ and $g_i:\bar{m_i}\to\ov{j_i}$, $1\leq i\leq k$, 
be nondegenerate.
For each diagram $D$ of type $(f,m_1,\ldots,m_k)$, let $h_D$ be the composite
\[
\overline{m-k+\scriptstyle\sum m_i}\labarrow{b}\textstyle\coprod\overline{m_i}
\labarrow{\coprod g_i}
\textstyle\coprod \overline{j_i}
\labarrow{\omega}
\overline{\scriptstyle\sum j_i}
\]
where $\omega$ takes $q\in\overline{j_i}$ to $j_1+\cdots+j_{i-1}+q$.
Let $\eta(D)$ denote
\[
\sum_{1\leq i<i'\leq k} (m_i-j_i)||f^{-1}(i')||
+
\sum_{i\in\bar{k}} \epsilon(g_i,\A^i) 
\]
\end{definition}

Finally, we have

\begin{proposition}
\label{gamma}
\[
\la f\ra (\la g_1\ra ,\ldots,\la g_k\ra )=
\sum_D (-1)^{\eta(D)}\,\la h_D\ra 
\]
where the sum is taken over all diagrams $D$ of type $(f,m_1,\ldots,m_k)$.
\end{proposition}

A conceptual way to remember the sign is that (in the terminology of Remark 
\ref{Fri1}(f)) it comes from permuting the $*$'s corresponding to $f$ and the 
$g_i$ into the order corresponding to $h$.

We defer the proof of Proposition \ref{gamma} to Section \ref{sec6}.

It only remains to prove part (c) of Theorem \ref{revS2}.  
Consider the chain homotopy
\[
s:\S(k)_q\to\S(k)_{q+1}
\]
which places a 1 at the beginning of each sequence; if the sequence already
begins with a 1 then the new sequence is degenerate so $s$ takes it to zero.
This chain homotopy has the property that
\[
\partial s+s \partial =\text{id} +\iota r
\]
where $\iota:\S(k-1)\to\S(k)$ places a 1 at the beginning of each sequence and
increases each of the original entries by 1, and $r:\S(k)\to\S(k-1)$ takes a
sequence to zero unless it begins with a 1 and has no other 1's, in
which case $r$ removes the 1 and decreases each of the remaining entries by 1.
Since $r$ is an epimorphism and $\iota$ is a monomorphism this implies that
$\S(k)$ has the same homology as $\S(k-1)$, so the desired result follows by
induction.

\begin{remark}
\label{revM4}
The chain homotopy that was used in this proof is due to Benson \cite[page 
147]{Benson}.
\end{remark}

\section{The chain operads $\S_n$}
\label{sec3}

The normalized singular chain functor
$S_*$ takes topological operads to chain operads (\cite[page 25]{KM}).
Let $\C_n$ be the little $n$-cubes operad (\cite{MayG},\cite{BV}).
In this section
we define a suboperad $\S_n$ of $\S$ which will turn out to be
quasi-isomorphic (in the category of chain operads over ${\mathbb Z}$) to 
$S_*\C_n$.

\begin{definition} 
\label{complexity}
Let $T$ be a finite totally ordered set, let $k\geq 2$, and let
$f:T\rightarrow \bar{k}$.
We define the {\it complexity} of $f$ as follows.  
If $k$ is $0$ or $1$ the complexity is 0.
If $k=2$ let $\sim$ be
the equivalence relation on $T$ generated by
\[
a\sim b \mbox{ if $a$ is adjacent to $b$ and $f(a)=f(b)$}
\]
and define the complexity of $f$ to be the number of equivalence classes
minus 1.  If $k> 2$  define the complexity of $f$ to be the maximum of the
complexities of the restrictions $f|_{f^{-1}(A)}$ as $A$ ranges over the
two-element subsets of $\bar{k}$.
\end{definition}

The motivation for this definition is given in  Remark \ref{motivation}.

\begin{definition}
For each $n\geq 1$ and $k\geq 0$, let $\S_n(k)$ be the 
sub-graded-abelian-group of $\S(k)$
generated by the nondegenerate $f$ with complexity $\leq n$.  Let $\S_n$ 
denote the collection $\S_n(k)$, $k\geq 0$.
\end{definition}

\begin{proposition}
\label{sub}
$\S_n$ is a sub-chain-operad of $\S$.
\end{proposition}

For the proof we need a lemma whose proof is left to the reader.

\begin{lemma}
\label{cx}
Given a commutative diagram
\[
\xymatrix{
& T'
\ar[dl]_{f'}
\ar[dd]^{g}
\\
\bar{k}
& \\
& T
\ar [ul]^{f}
}
\]
with $g$ ordered,
the complexity of $f'$ is $\leq$ the complexity of $f$.  If $g$ is a
surjection then the complexity of $f'$ is equal to the complexity of $f$.
\end{lemma}

\begin{proof}[Proof of Proposition \ref{sub}.]
It's easy to see that the action of 
$\Sigma_k$ on $\S_n(k)$ preserves complexity.  The fact that $\partial$
preserves complexity is immediate from Lemma \ref{cx}. Now suppose that 
$f:\bar{m}\rightarrow\bar{k}$ 
and $g_i:\overline{m_i}\rightarrow\overline{j_i}$ have complexity $\leq n$
and choose a diagram $D$ of type $(f,m_1,\ldots,m_k)$:
\[
\xymatrix{
\overline{m-k+\sum m_i}
\ar[r]^-{b}
\ar[d]_{a}
&
\coprod\overline{m_i}
\ar[d]_{\chi}
\\
\bar{m}
\ar[r]^{f}
&
\bar{k}
}
\]
It suffices to show that $h_D$ has complexity $\leq n$, and
for this we need to 
show that the restriction of $h_D$ to $h_D^{-1}(A)$ has complexity $\leq n$
whenever $A$ is a two-element subset of 
$\overline{\textstyle\sum j_i}$.  
There are two cases: either $\omega^{-1}(A)$ is
contained in some $\overline{j_i}$ or not.  
If $\omega^{-1}(A)\subset
\overline{j_i}$ for some $i$ then the complexity of $h_D|_{h_D^{-1}(A)}$ is 
equal to the complexity of $(\omega\circ g_i\circ b)|_{h_D^{-1}(A)}$, and 
this is less than or equal to the complexity of $g_i$ by Lemma \ref{cx} 
(since $b$ is order-preserving on $b^{-1}\overline{m_i}$).

If $\omega^{-1}(A)$ is not contained in any $\overline{j_i}$ then the 
complexity of $h_D|_{h_D^{-1}(A)}$ is the same as that of 
$\chi\circ b|_{h_D^{-1}(A)}$.
But $\chi\circ b=f\circ a$, and the complexity of $f\circ a|_{h_D^{-1}(A)}$
is less than or equal to the complexity of $f$ by Lemma \ref{cx} (since $a$ is
order-preserving).
\end{proof}

\begin{theorem}
\label{big}
$\S_n$ is quasi-isomorphic, in the category of chain operads, to
$S_*\C_n$. 
\end{theorem}

This means that there is a sequence 
\[
\S_n \leftarrow \cdots \rightarrow S_*\C_n
\]
of chain-operad maps which are quasi-isomorphisms.

Theorem \ref{big} will be proved in Section \ref{sec7}.

\begin{remark}
\label{motivation}
(Motivation for the definition of complexity.)
The definition of complexity when $k=0$ is motivated by the fact that $\C_n(0)$
is a point for all $n$.  The definition when $k=1$ is motivated by the
requirement that the unit of $\S$ should be in $\S_n$.  The definition when
$k=2$ is motivated by the fact that $\C_n(2)$ is homotopic to $S^{n-1}$; note 
that $\S_n(2)$ is isomorphic to the cellular chain complex of 
the usual ${\mathbb Z}/2$-equivariant CW structure on $S^{n-1}$.
The definition of complexity when $k>2$ is not as easy to motivate, and it is
possible that there are other ways of choosing the complexity filtration so
that the analog of Theorem \ref{big} is true (note, however, how easy Lemma
\ref{cx} is to prove for the filtration we have given).  Our choice was 
suggested by 
the filtration used in \cite{Smith} (the original motivation for the filtration
in \cite{Smith} is that it is the simplest filtration
compatible with the ``degeneracy maps'' from the $k$-th space of the operad to
the second space of the operad).
\end{remark}

\section{Deligne's Hochschild cohomology conjecture.}

\label{sec3a}

In this section we describe a natural action of $\S_2$ on the normalized 
Hochschild cochain complex.  Intuitively, the reason there is such an action 
is that (i) the cup product and brace operations (see \cite{Getzler}, 
\cite{Kad}) 
on Hochschild cochains satisfy the same relations as the sequence 
operations 12 and $12131\cdots 1$, and (ii) these sequence operations 
generate $\S(2)$ as an operad.

Let $R$ be an associative ring. 
Recall that the $p$-th Hochschild cochain group consists of the homomorphisms
of abelian groups
\[
x:R^{\otimes p}\to R
\]
where $R^{\otimes 0}$ is interpreted as $\mathbb Z$.
The normalized cochain group $C^p (R)$ consists of the Hochschild cochains 
whose composites with each of the maps
\[
R^{\otimes i}\otimes {\mathbb Z} \otimes R^{\otimes (p-i-1)}
\labarrow{1\otimes \iota \otimes 1}
R^{\otimes p}
\]
are zero, where $\iota$ is the unit of $R$.
The differential $C^p(R)\to C^{p+1}(R)$ is defined by
\begin{multline}
\notag
(\partial x)(r_1\otimes\cdots\otimes r_{p+1})
=
r_1\cdot x(r_2\otimes\cdots\otimes r_{p+1}) \\
+
\sum_{0<i<p+1} x(r_1\otimes\cdots\otimes r_i r_{i+1}\otimes\cdots\otimes
r_{p+1})
+
(-1)^{p+1}x(r_1\otimes\cdots\otimes r_p)\cdot r_{p+1}
\end{multline}
There is also a cup product in $C^*(R)$: if $x\in C^p(R)$ and $y\in C^q(R)$
then 
\[
(x\smallsmile y)(r_1\otimes\cdots\otimes r_{p+q})
=
x(r_1\otimes\cdots\otimes r_p)\cdot y(r_{p+1}\otimes\cdots\otimes r_{p+q})
\]
Finally, if $x\in
C^{k}(R)$ and $x_i\in C^{j_i}(R)$ for $1\leq i\leq k$ we define
$x(x_1,\ldots,x_k)\in C^{j_1+\cdots+j_k}(R)$ to be the composite
\[
R^{j_1+\cdots+j_k}\labarrow{x_1\otimes\cdots\otimes x_k} R^{\otimes k}
\labarrow{x} R
\]

Now suppose given a sequence $\x=(x_1,\ldots,x_k)$
of normalized Hochschild cochains, 
a finite totally ordered set $T$, and a map 
$g:\bar{q}\to\bar{k}$ 
of complexity $\leq 2$
with the property that $|g^{-1}(i)|=|x_i|+1$ for each
$i$.  By a {\it segment} we mean a subset $S$ of $T$ such that $g$ has the same
value on the minimum and maximum elements of $S$, and by a {\it maximal segment}
we mean a segment which is not properly contained in any other segment.
Let $S_1,\ldots,S_r$ denote the maximal segments;
the fact that $g$ has complexity $\leq 2$ implies that the sets $S_j$ are
disjoint.
We define a normalized Hochschild cochain $g(\x)$ inductively as follows.

\begin{description}
\item[{\rm (i)}] If $T$ is empty then $g(\x)$ is the identity cochain in 
$C^1(R)$.

\item[{\rm (ii)}] If $T$ has a single element $t$ then $g(\x)$ is $x_{g(t)}$.

\item[{\rm (iii)}]  If $r>1$ then $g(\x)$ is the cup product 
$g|_{S_1}(\x)\smallsmile\cdots\smallsmile
g|_{S_r}(\x)$

\item[{\rm (iv)}] If $r=1$ and $|T|>1$ let $i$ be the 
value of $g$ at
the minimum and maximum values of $T$ and let $t_1,\ldots,t_{|x_i|+1}$ be the
elements of $g^{-1}(i)$ in increasing order.  For each $j< |x_i|+1$ let $A_j$
be the set $\{t_j < t < t_{j+1}\}$ (which may be empty) and let 
\[
g(\x)=x_i(g|_{A_1}(\x),\ldots,g|_{A_{|x_i|}}(\x))
\]

\end{description}

Now let $f:\bar{m}\to\bar{k}$ be nondegenerate, and let $\x=(x_1,\ldots,x_k)$
be a $k$-tuple in $C^*(R)$.  We define 
\[
\theta_k(f,x_1,\ldots,x_k)=
\sum_E (-1)^{\epsilon'(f,\A)} (f\circ\beta)(\x)
\]
where $E$ runs through the special diagrams
\[
\xymatrix{
\ov{\Sum |x_i| +k-m+1}
&
\ov{k+\Sum|x_i|}
\ar[l]_-{\alpha}
\ar[r]^-{\beta}
&
\bar{m}
}
\]
and $\A$ denotes the overlapping partition of $\ov{k-m+\Sum |x_i|}$ with
$\A_j=\alpha(\beta^{-1}(j))$; the sign is given by
\[
\epsilon'(f,\A)=\epsilon(f,\A)+(m-k)\sum|x_i|+|\{(j,j'):f(j)=f(j')\}|
\]
The term corresponding to $E$ is counted as zero unless, for each $i$,
$\alpha|_{\beta^{-1}f^{-1}(i)}$ is a monomorphism and
$|\beta^{-1}f^{-1}(i)|=|x_i|+1$.

\begin{theorem}
\label{De}
The maps $\theta_k$ give $C^*(R)$ a natural structure of algebra over $\S_2$.
\end{theorem}

The proof is similar to the proofs of Propositions \ref{diff}, \ref{perm} and 
\ref{gamma} and is left to the reader.

\section{Proof of Theorem \ref{big}}

\label{sec7}

\begin{definition}
For each $k\ge0$, 
let $P_{2}\overline{k}$ be the set of subsets of
$\overline{k}$ that have two elements. 
\end{definition}

\begin{definition}  
Let $\mathcal I(k)$ be the set whose elements are pairs
  $(b,T)$, where $b$ is a function from $P_{2}\overline{k}$ to the nonnegative
integers and $T$ is a total order of $\ov{k}$. We give 
$\mathcal I(k)$ the partial order for which
  $(a,S)\le (b,T)$ 
  if $a(\{i,j\})\le b(\{i,j\})$ for each $\{i,j\}\in
  P_{2}\overline{k}$ and $a(\{i,j\})< b(\{i,j\})$ for each $\{i,j\}$
  with $i<j$ in the order $S$ but $i>j$ in the order $T$. Let
  $\I_n(k)$ be the subset of pairs $(b,T)$ 
  such that $b\{i,j\}<n$ for each $\{i,j\}\in P_{2}\ov {k}$. The
  set $\I_n(k)$ inherits an order from $\I(k)$.
\end{definition}

The sequence $\I(k)$ is an operad in the category of
partially ordered sets with the following structure maps \cite{Berger}. The 
right action of
$\Sigma_{k}$ on $\I(k)$ is given by 
\[
(b,T)\rho=(b\circ\rho_{2},T\rho)
\]
where $\rho_{2}\colon P_{2}\ov k \to P_{2}\ov k$ is the function
$\rho_{2}(\{i,j\})=\{\rho(i),\rho(j)\}$ and where $i<j$ in the
total order $T\rho$ if $\rho(i)<\rho(j)$ in the total order $T$.
The operad composition
\[
\I(k)\times\I(a_{1})\times\dots\times\I(a_{k})\to \I(\Sigma a_{i})
\] 
sends the $(k+1)$-tuple of pairs
$((b,T);(b_1,T_{1}),\dots,(b_k,T_{k}))$ with $(b,T)\in\I(k)$ and
$(b_{i},T_{i})\in\I(a_{i})$ to the pair
$(b(b_1,\dots,b_{k}),T(T_1,\dots,T_k))$ in $\I(a)$ with $a=\Sigma
a_{i}$, where the value of $b(b_1,\dots,b_k)$ at $\{r,s\}$ is
$b_{i}(\{r,s\})$ if $\{r,s\}\subset\ov a_{i}$ and the value at
$\{r,s\}$ is $b(\{i,j\})$ if $r\in\ov a_{i}$, $s\in\ov a_{j}$ and
$i\neq j$, and where $T(T_1,\dots,T_k)$ is the total order of $\ov
a=\coprod \ov a_{i}$ for which $r<s$ if $\{r,s\}\subset\ov a_{i}$ and
$r<s$ in the order $T_{i}$ on $\ov a_{i}$ or if $r\in\ov
a_{i}$, $s\in\ov a_{j}$ and $i<j$. For each $n$ the sequence of
partial orders $\I_n(k)$ is a suboperad of the sequence  $\I(k)$.

Given a category $\I$ let us write $\Ne\I$ for the nerve of $\I$
and $C_* \I$ for the normalized chains of the simplicial set $\Ne\I$.

\begin{proposition}
  The chain operad $C_{*}\I_n(k)$ is quasi-isomorphic to the
  chain operad $S_{*}\C_{n}$.
\end{proposition}

\begin{proof}
  Berger has shown \cite{Berger} that the operad of spaces obtained as the 
geometric
  realization of the nerves $|\Ne\mathcal I_n(k)|$ is weakly
  equivalent to the operad of little $n$-cubes. Now apply the normalized 
singular
  chains functor to get a quasi-isomorphism of chain operads between
  $S_{*}|\Ne\mathcal I_n(k)|$ and $S_{*}\C_{n}$. The natural map
  $C_{*}\I_n(k)\to S_{*}|\Ne\mathcal I_n(k)|$ is a 
  quasi-isomorphism of chain operads and the proof is complete. 
\end{proof}

Next, a technique of Berger \cite{Berger} can be used to prove that the 
operad $\S_{n}$ is quasi-isomorphic to the operad $C_{*}\I_{n}$.

\begin{definition}
  Let $f\colon \ov{m}\to \ov k$ be a surjective map. Then
  $(b_f,T_{f})\in\I(k)$ is the pair where $b_{f}(\{i,j\})$ is one
  less than the complexity of the restriction of $f$ to a map
  $f^{-1}(\{i,j\})\to\{i,j\}$ and where $i<j$ in the total order
  $T_{f}$ if the smallest element of $f^{-1}(i)$ is less than the
  smallest element of $f^{-1}(j)$.
\end{definition}

\begin{definition}
  For $(b,T)\in\I(k)$, the chain complex $\S(b,T)$ is the
  subcomplex of $\S(k)$ generated by those surjective
  functions $f\colon\ov{m}\to\ov k$ with $(b_f,T_{f})\le
  (b,T)$. 
\end{definition} 

\begin{proposition}
  The chain complex $\S(b,T)$ is contractible.
\end{proposition}

\begin{proof}
  For $i\in\ov k$ let $i_{*}\S(k-1)$ be the subcomplex of $\S(k)$ 
  generated by the sequences that begin with $i$ and have no other
  occurrences of $i$.  In the proof of Theorem \ref{revS2}(c) we constructed a 
chain
  homotopy $s\colon\S(k)\to\S(k)$ that gives a deformation
  retraction of $\S(k)$ onto $1_{*}\S(k-1)$.  For $i\neq1$ let
  $\rho_{i}$ be the transposition $(1,i)$ and let $s_{i}=s\circ\rho_{i}$. Then 
$s_{i}$ is
  a chain homotopy that gives 
  a deformation retraction of $\S(k)$ onto $i_{*}\S(k)$. If a subcomplex $C$
  of $\S(k)$ is invariant under $s_i$ then there is a deformation
  retraction of $C$ onto $C\cap i_{*}\S(k-1)$.
  In particular, if $i$ is the first element of the total order $T$,
  then $\S(b,T)$ is invariant under $s_i$ and its deformation
  retract is isomorphic to $\S(\delta_i^{*}b,\delta_i^{*}T)$ where
  $\delta_{i}\colon\ov{k-1}\to\ov k$ is the unique order-preserving
  monomorphism that does not have $i$ in its image, $\delta_{i}^{*}b$
  is the restriction of $b$ to $P_{2}\ov {k-1}$ and $\delta_{i}^{*}T$
  is the pullback of $T$ to a total order of $\ov {k-1}$. It then
  follows by induction that $\S(b,T)$ is contractible.
\end{proof}

Let ${\mathbf C}$ denote the category of
non-negatively graded chain complexes of abelian groups. 

\begin{definition}
Let $D_n(k)\colon\I_n(k)\to{\mathbf C}$ be the diagram of
  chain complexes given by 
\[
D_n(k)(b,T)=\S(b,T)
\]
for $(b,T)\in\I_n(k)$.
\end{definition}

A function $f\colon \ov{m}\to \ov k$ is in $\S(b,T)$ if and only
$(b_{f},T_{f})\le (b,T)$. It follows that the natural map
\[
\colim_{\I_n(k)} D_n(k)\to \S_n(k)
\]
is an isomorphism. We wish to study the homotopy colimit of $D_n(k)$
and so we recall (from \cite[Section 20.1]{HH}) the definition of the 
homotopy colimit in the category of chain complexes.  
%Let ${\mathbf C}$ be the category of
%non-negatively graded chain complexes of abelian groups.  
%For a
%category $\I$, the nerve $\Ne\I$ is a simplicial set and let $C_{*}I$
%be the normalized chain complex of $\Ne\I$. 
For a partially ordered set
$\I$, let $a/\I$ denote the suborder of all elements of $b\in \I$ such
that $a\ge b$.  Notice that $C_{*}(-/\I)$ is a contravariant functor
from $\I$ to ${\mathbf C}$

\begin{definition}
  Let $\I$ be a partially ordered set and let $D\colon \I\to \C$ be a
  diagram of chain complexes.  The \emph{homotopy colimit} of $D$ is
  the coequalizer in ${\mathbf C}$
\[
\bigoplus_{a<b}C_{*}(b/\I)\otimes Da
\stackrel{\textstyle{\to}}{\to}\bigoplus_{a}C_{*}(a/\I)\otimes Da \to \hocolim D
\]
\end{definition}

As a special case, if $D$ is the constant functor with value $\mathbb
Z$ then $\hocolim D=C_{*} \I$. 

\begin{proposition}\label{prop-cofibrant-D}
  Let $\I$ be a partially ordered set and let $D\colon \I\to \C$ be a
  diagram of chain complexes. If the map $\colim_{a<b}Da \to Db$ is a
  monomorphism for each $b\in \I$ then the natural map $\hocolim D\to
  \colim D$ is a quasi-isomorphism.
\end{proposition}

\begin{proof}
This is a general property of homotopy colimits, see 
\cite[Theorem 20.9.1]{HH}
Give $\C$, the category of chain complexes, the model structure for 
which the cofibrations are the monomorphisms and the weak 
equivalences are the quasi-isomorphisms.  Then $D\colon \I\to \C$ is a 
Reedy cofibrant diagram and $\hocolim D\to \colim D$ is a 
quasi-isomorphism.
\end{proof}

For $\rho\in\Sigma_{k}$, the action of $\rho$ on $\S(k)$ restricts
to a natural map 
\[
\rho^{*}\colon\S(b,T)\to\S((b,T)\rho). 
\]
This is a natural transformation of functors on $\I_n(k)$ and
induces a map of homotopy colimits
\[
\rho^*:\hocolim D_n(k)\to\hocolim D_n(k)
\]
Also, the operad composition of $\S$ induces a natural transformation
\[
\S(b,T)\otimes \bigotimes_{i}\S(b_{i},T_{i})\to
\S(b(b_1,\dots,b_k),T(T_1,\dots,T_{k}))
\]
and this induces a map
\[
\hocolim D_n(k)\otimes\bigotimes_{i}\hocolim D_n(a_{i})\to
\hocolim D_n(\Sigma a_{i})
\]
Chasing diagrams gives:

\begin{proposition}
  For each $n\ge1$ the sequence of chain complexes $\hocolim D_n(k)$
  is an operad.
\end{proposition}

Now we can complete the proof of Theorem \ref{big}.
We have a diagram of chain operads
\[
C_{*}\I_n(k)=\hocolim_{\I_n(k)}\mathbb Z\leftarrow
\hocolim_{\I_n(k)} \S(b,T) \to \colim_{\I_n(k)}
\S(b,T)=\S_n(k).
\]
Since the chain complexes $\S(b,T)$ are contractible, the left hand
map is a quasi-isomorphism of chain operads. Since the diagrams $D_n(k)$
satisfy the condition in proposition \ref{prop-cofibrant-D} the right
hand map is also a quasi-isomorphism of chain operads.

\section{Proof of Proposition \ref{perm}.}
\label{sec4}

Recall the homomorphism 
\[
[f]:\S_* X\to (S_* X)^{\otimes k}
\]
from Definition \ref{la}(a).  To prove Proposition \ref{perm} it suffices to
show
\begin{equation}
\label{revT1}
\rho\circ [f]=(-1)^{\zeta(f,\rho)}[ \rho\circ f]
\end{equation}

Recall that when $A$ is a set we are using $||A||$ to denote the cardinality of
$A$ minus 1.

Let us fix an $i$ with $1\leq i\leq k$.
It suffices to verify (\ref{revT1}) when $\rho$ is the 
transposition permuting 
$i$ and $i+1$.  We need to show that
\begin{equation}
\label{1}
\epsilon(f,\A)+||\tcoprod_{f(j)=i} \, A_j||\,||\tcoprod_{f(j)=i+1} \, A_j||
\equiv \zeta(f,\rho) +\epsilon(\rho\circ f,\A) \qquad\mbox{mod 2}
\end{equation}
Unwinding the definitions, the left-hand side of (\ref{1}) becomes
\begin{multline*}
\underset{f(j)>f(j')}{\underset{j<j'}\sum} ||A_j||\,||A_{j'}||+
\sum_{j=1}^m ||A_j||(\tau_f(j)-f(j)) +  \\
\biggl(
||f^{-1}(i)||
+
\sum_{j\in f^{-1}(i)} ||A_j||
\biggr)
\biggl(
||f^{-1}(i+1)||
+
\sum_{j\in f^{-1}(i+1)} ||A_j||
\biggr)
\end{multline*}
and the right-hand side becomes 
\[
||f^{-1}(i)||\,||f^{-1}(i+1)|| + 
\underset{\rho f(j)>\rho f(j')}{\underset{j<j'}\sum} 
||A_j||\,||A_{j'}||+
\sum_{j=1}^m ||A_j||(\tau_{\rho f}(j)-\rho f(j))
\]
The fact that the left- and right-hand sides of (\ref{1}) are congruent mod 2 
follows easily from the fact that
\[
\tau_{\rho \circ f}(j)-(\rho \circ f)(j)=
\left\{
\begin{array}{ll}
\tau_f(j)+||f^{-1}(i+1)|| & \mbox{if $f(j)=i$} \\
\tau_f(j)-||f^{-1}(i)|| & \mbox{if $f(j)=i+1$} \\
\tau_f(j)-f(j) & \mbox{otherwise}
\end{array}
\right.
\]

\section{Proof of Proposition \ref{diff}.}
\label{sec5}

For this proof it will be convenient to use special diagrams instead of
overlapping partitions (see Definition 
\ref{special} and Lemma \ref{part}) so we begin by reformulating Definitions
\ref{revT3} and \ref{la} in this language.

\begin{definition}
\label{eps}
Given a special diagram 
\[
E: \quad
\xymatrix{
\{0,\ldots,p\}
&
\overline{p+m}
\ar[l]_-{\alpha}
\ar[r]^-{\beta}
&
\bar{m}
}
\]
and a nondegenerate $f:\bar{m}\to\bar{k}$,
let 
\[
\epsilon(f,E)=
\underset{f(j)>f(j')}{\underset{j<j'}\sum}
||\beta^{-1}(j)||\,||\beta^{-1}(j')||
+
\sum_{j=1}^m ||\beta^{-1}(j)||(\tau_f(j)-f(j))
\]
\end{definition}

Using Lemma \ref{part} 
we may rewrite Definition \ref{la} as
\begin{equation}
\label{ss2}
\sigma[f]=\sum_E (-1)^{\epsilon(f,E)} 
\bigotimes_{i=1}^k \sigma(\alpha(\beta^{-1}f^{-1}(i)))
\end{equation}
where the sum is taken over all special diagrams 
\[
E: \quad
\xymatrix{
\{0,\ldots,p\}
&
\overline{p+m}
\ar[l]_-{\alpha}
\ar[r]^-{\beta}
&
\bar{m}
}
\]
and the term 
corresponding to $E$ is counted as zero unless $\alpha$ is a
monomorphism on $\beta^{-1}f^{-1}(i)$ for each $i$.

Now fix $f:\bar{m}\to\bar{k}$ and $\sigma:\Delta^p\to X$.  To prove Proposition
\ref{diff} it suffices to show
\begin{equation}
\label{3}
\partial((\sigma)[f])=
(\partial\sigma)[f]
+(-1)^p
\sum_{j=1}^m (-1)^{\tau_f(j)-f(j)} \sigma[f_{\bar{m}-\{j\}}]
\end{equation}

First we analyze the left side of equation (\ref{3}).  By equation 
(\ref{ss2}) it
is equal to
\begin{equation}
\label{4}
\partial\Bigl(
\sum_E (-1)^{\epsilon(f,E)}
\bigotimes_{i=1}^k \sigma(\alpha(\beta^{-1}f^{-1}(i)))
\Bigr)
\end{equation}
For each pair $(E,q)$ with $q\in\overline{p+m}$ let us write $E_q$
for the diagram
\[
\xymatrix{
\alpha(\overline{p+m}-\{q\})
&
\overline{p+m}-\{q\}
\ar[l]_-{\alpha_q}
\ar[r]^-{\beta_q}
&
\beta(\overline{p+m}-\{q\})
}
\]
where $\alpha_q$ and $\beta_q$ are the restrictions of $\alpha$ and $\beta$ (we
are not claiming that all of the diagrams $E_q$ obtained in this way are
special).  Let $f_q$ be the restriction of $f$ to
$\beta(\overline{p+m}-\{q\})$.
Then expression (\ref{4}) can be rewritten as
\begin{equation}
\label{4a}
\sum_{E,q} (-1)^{\epsilon(f,E)}
\bigotimes_{i=1}^k 
(-1)^{\kappa(f,E,q)}
\sigma(\alpha_q(\beta_q^{-1}f_q^{-1}(i)))
\end{equation}
where 
\begin{equation}
\label{5}
\kappa(f,E,q)
=
\sum_{i<f(\beta(q))} ||\beta^{-1}f^{-1}(i)||
\ \ +\ 
|\{q'<q : f(\beta(q'))=f(\beta(q))\}|
\end{equation}
Next let us divide the pairs $(E,q)$ into four types:

\begin{itemize}
\item $(E,q)$ is of type I if $1<q<p+m$ and 
$\alpha(q-1)<\alpha(q)=\alpha(q+1)$.
\item $(E,q)$ is of type II if $1<q<p+m$ and 
$\alpha(q-1)=\alpha(q)<\alpha(q+1)$.
\item $(E,q)$ is of type III if $1<q<p+m$ and 
$\alpha(q-1)<\alpha(q)<\alpha(q+1)$, or if $q=1$ and $\alpha(1)<\alpha(2)$, or
if $q=p+m$ and $\alpha(p+m-1)<\alpha(p+m)$.
\item $(E,q)$ is of type IV if $1<q<p+m$ and 
$\alpha(q-1)=\alpha(q)=\alpha(q+1)$, or if $q=1$ and $\alpha(1)=\alpha(2)$, or
if $q=p+m$ and $\alpha(p+m-1)=\alpha(p+m)$.
\end{itemize}

To complete the proof of Proposition \ref{diff} it suffices to show

\begin{lemma}
\label{atlast}
In expression (\ref{4a}),

(a) the terms of type I cancel the terms of type II

(b) the terms of type III add up to 
$(\partial\sigma)[f]$

(c) the terms of type IV add up to 
$(-1)^p\sum_{j=1}^m (-1)^{\tau_f(j)-f(j)} \sigma[f_{\bar{m}-\{j\}}]$
\end{lemma}

\begin{proof}[Proof of Lemma \ref{atlast}.]
First we simplify equation (\ref{5}). We have
\begin{equation}
\notag
\begin{split}
\kappa(f,E,q) =\,&q-1-
\underset{f(j)\neq f(\beta(q))}{\underset{j<\beta(q)}\sum} |\beta^{-1}(j)|
+
\underset{f(j)<f(\beta(q))}\sum |\beta^{-1}(j)| \ \ - (f(\beta(q))-1)
\\
=\,&q-f(\beta(q))-
\underset{f(j)> f(\beta(q))}{\underset{j<\beta(q)}\sum} |\beta^{-1}(j)|
+
\underset{f(j)< f(\beta(q))}{\underset{j>\beta(q)}\sum} |\beta^{-1}(j)|
\\
=\,&q-f(\beta(q))-
\underset{f(j)> f(\beta(q))}{\underset{j<\beta(q)}\sum} ||\beta^{-1}(j)||
+
\underset{f(j)< f(\beta(q))}{\underset{j>\beta(q)}\sum} ||\beta^{-1}(j)||
\\
\,&-|\{j<\beta(q) : f(j)>f(\beta(q))\}|+|\{j>\beta(q) : f(j)<f(\beta(q))\}|
\end{split}
\end{equation}
Using the fact that 
\[
\tau_f(\beta(q))=|\{j\leq\beta(q) : f(j)\leq f(\beta(q))\}|
+ |\{j>\beta(q) : f(j)<f(\beta(q))\}|
\]
and the fact that
\begin{equation}
\label{5z}
\alpha(q)=q-\beta(q)
\end{equation}
the last expression simplifies to
\begin{equation}
\label{5a}
\kappa(f,E,q)
=
\alpha(q)+ \tau_f(\beta(q))-f(\beta(q))
-
\underset{f(j)>f(\beta(q))}{\underset{j<\beta(q)}\sum} ||\beta^{-1}(j)||
+
\underset{f(j)<f(\beta(q))}{\underset{j>\beta(q)}\sum} ||\beta^{-1}(j)||
\end{equation}

From equation (\ref{5a}) it is easy to see that if $(E,q)$ is of type I, II 
or III we have
\begin{equation}
\label{6}
\epsilon(f,E)+\kappa(f,E,q)\equiv \alpha(q)+\epsilon(f,E_q)\quad\mbox{mod 2}
\end{equation}
Now suppose $E$ is of type I.  Define
\[
E':\quad \xymatrix{
\{0,\ldots,p\}
&
\overline{p+m}
\ar[l]_-{\alpha'}
\ar[r]^-{\beta'}
&
\bar{m}
}
\]
by letting 
\[
\alpha'(r)=
\left\{
\begin{array}{ll}
\alpha(r)&\mbox{if $r\neq q$} \\
\alpha(q)-1&\mbox{if $r=q$}
\end{array}
\right.
\]
and
\[
\beta'(r)=
\left\{
\begin{array}{ll}
\beta(r)&\mbox{if $r\neq q$} \\
\beta(q)+1&\mbox{if $r=q$}
\end{array}
\right.
\]
Then $E'$ is of type II, and $E'_q$ is the same diagram as $E_q$, but
$\alpha'(q)=\alpha(q)-1$.  Now equation (\ref{6}) implies that the terms in the
expression (\ref{4a}) corresponding to $E$ and $E'$ cancel, which completes the
proof of part (a).

For part (b), let $E$ be a diagram of type III.  In this case 
\[
\alpha(\overline{p+m}-\{q\})=\{0,\ldots,p\}-\{\alpha(q)\},
\]
so we want to show
\begin{multline*}
p(m-k)+\epsilon(f,E)+\kappa(f,E,q)
\equiv \\
(m-k)+(p-1)(m-k)+\alpha(q)+\epsilon(f,E_q) \quad\mbox{mod 2}
\end{multline*}
which is immediate from equation \eqref{6}

For part (c), Let $E$ be a diagram of type IV.  
In this case 
\[
\beta(\overline{p+m}-\{q\})=\bar{m}-\{\beta(q)\}
\]
so we need to show
\begin{equation}
\label{7}
\epsilon(f,E)+\kappa(f,E,q)\equiv
p+\tau_f(\beta(q))-f(\beta(q))+\epsilon(f_q,E_q)
\quad\mbox{mod 2}
\end{equation}
Expanding the definitions and making the obvious
cancellations, equation (\ref{7}) reduces to
\begin{multline}
\label{8}
\alpha(q)-
\underset{f(j)>f(\beta(q))}{\underset{j<\beta(q)}\sum} ||\beta^{-1}(j)||
+
\underset{f(j)<f(\beta(q))}{\underset{j>\beta(q)}\sum} ||\beta^{-1}(j)||
\ \ +p+
\sum_{j=1}^m||\beta^{-1}(j)||(\tau_f(j)-f(j))
\\
\equiv
\sum_{j\neq q}||\beta^{-1}(j)||(\tau'_f(j)-f(j))
\end{multline}
Here 
\[
\tau'_f(j)=
\left\{
\begin{array}{ll}
\tau_f(j)-1 &\mbox{if $j>\beta(q)$ and $f(j)=f(\beta(q))$, or if
$f(j)>f(\beta(q))$} \\
\tau_f(j) &\mbox{otherwise}
\end{array}
\right.
\]
so equation (\ref{8}) reduces to 
\[
\alpha(q)+p\equiv 
\underset{j>\beta(q)}\sum ||\beta^{-1}(j)|| \quad\mbox{mod 2}
\]
and this follows from equation (\ref{5z}) and the fact that
\[
\underset{j>\beta(q)}\sum ||\beta^{-1}(j)||=
\underset{j>\beta(q)}\sum |\beta^{-1}(j)| \ \ -(m-\beta(q))
=
p+m-q-(m-\beta(q))
\]
\end{proof}

\section{Proof of Proposition \ref{gamma}.}
\label{sec6}

Fix nondegenerate maps $f:\bar{m}\to\bar{k}$ and $g_i:\bar{m_i}\to\bar{j_i}$ 
for $1\leq i\leq k$. Also fix $\sigma:\Delta ^p\to X$.  We need to show
\begin{multline}
\label{s1}
\la f\ra (\la g_1\ra ,\ldots,\la g_k\ra )(x_1\otimes\ldots\otimes
x_{j_1+\cdots+j_k})(\sigma)
= \\
\sum_D (-1)^{\eta(D)}\,\la h_D\ra
(x_1\otimes\ldots\otimes
x_{j_1+\cdots+j_k})(\sigma)
\end{multline}
where $D$ runs through the diagrams of type $(f,m_1,\ldots,m_k)$.

For the proof we will work with special diagrams instead of 
overlapping partitions (see Definition \ref{special} and Lemma \ref{part}).  

First we apply equation \eqref{ss2} to expand the right-hand side of 
\eqref{s1}:
\begin{multline}
\label{revTh1}
\sum_D (-1)^{\eta(D)}\,\la h_D\ra
(x_1\otimes\ldots\otimes
x_{j_1+\cdots+j_k})(\sigma)
\\
=
(x_1\otimes\cdots\otimes x_{j_1+\cdots+j_k})\Bigl(
\sum_D (-1)^{\eta(D)}\sum_F (-1)^{|h_D|+\epsilon(h,F)}
\bigotimes_{q=1}^{\sum j_i} \sigma(\phi(\psi^{-1}h_D^{-1}(q)))
\Bigr)
\end{multline}
where $F$ runs through the special diagrams of the form
\[
\xymatrix{
\{0,\ldots,p\}
&
\overline{p+m-k+\sum m_i}
\ar[l]_-{\phi}
\ar[r]^-{\psi}
&
\overline{m-k+\sum m_i}
}
\]

Next we apply equation
\eqref{ss2} twice to expand the left-hand
side of \eqref{s1}:
\begin{align}
\label{s2}
&\la f\ra (\la g_1\ra ,\ldots,\la g_k\ra )(x_1\otimes\ldots\otimes
x_{j_1+\cdots+j_k})(\sigma)
\\
\notag
&=
(x_1\otimes\cdots\otimes x_{j_1+\cdots+j_k})\Bigl(
(-1)^{|f|+\sum|g_i|}
\sum_E
(-1)^{\epsilon(f,E)}
\bigotimes_{i=1}^k
\\
\notag
&\quad\quad\quad\quad\quad
\Bigl(
\sum_{E_i}
(-1)^{
\sum_{i>i'} (m_{i'}-j_{i'})||\beta^{-1}f^{-1}(i)||
+\epsilon(g_i,E_i)
}
\bigotimes_{r=1}^{j_i}
\sigma(\alpha_i(\beta_i^{-1}g_i^{-1}(r)))
\Bigr)\Bigr)
\end{align}

Here $E$ runs through the special diagrams
\begin{equation}
\label{E}
\xymatrix{
\{0,\ldots,p\}
&
\overline{p+m}
\ar[l]_-{\alpha}
\ar[r]^-{\beta}
&
\bar{m}
}
\end{equation}
which satisfy
\begin{equation}
\label{m2}
\mbox{$\alpha$ is a monomorphism on $\beta^{-1}f^{-1}(i)$ for each 
$i\in\bar{k}$}
\end{equation}
(terms which do not satisfy this condition will be zero because they contain a
degenerate simplex)
and, for each choice of $E$ and $i$, $E_i$ runs through the special diagrams
\begin{equation}
\label{Ei}
\xymatrix{
\beta^{-1}f^{-1}(i)
&
\overline{||\beta^{-1}f^{-1}(i)||+m_i}
\ar[l]_-{\alpha_i}
\ar[r]^-{\beta_i}
&
\bar{m_i}
}
\end{equation}
which satisfy
\begin{equation}
\label{m3}
\mbox{$\alpha_i$ is a monomorphism on $\beta_i^{-1}g_i^{-1}(r)$ for each $r\in
\ov{j_i}$}
\end{equation}

In order to compare equations \eqref{revTh1} and $\eqref{s2}$, we must first
find the relationship between the indexing sets for the double sums.

Let ${\mathbf S}_1$ be the set of all $(k+1)$-tuples $(E,E_1,\ldots,E_k)$,
where $E$ is of the form \eqref{E} and satisfies \eqref{m2} and $E_i$ is of 
the form \eqref{Ei} and satisfies \eqref{m3} (${\mathbf S}_1$ is the indexing 
set for the double sum in equation \eqref{revTh1}).
Let 
${\mathbf S}_2$ be the set of all pairs $(D,F)$, where $D$ is a diagram 
\[
\xymatrix{
\overline{m-k+\sum m_i}
\ar[r]^-{b}
\ar[d]_{a}
&
\coprod\overline{m_i}
\ar[d]_{\chi}
\\
\bar{m}
\ar[r]^{f}
&
\bar{k}
}
\]
of type $(f,m_1,\ldots,m_k)$ and $F$ is a special diagram of the form
\[
\xymatrix{
\{0,\ldots,p\}
&
\overline{p+m-k+\sum m_i}
\ar[l]_-{\phi}
\ar[r]^-{\psi}
&
\overline{m-k+\sum m_i}
}
\]
satisfying
\begin{equation}
\label{m4}
\mbox{$\phi$ is a monomorphism on $\psi^{-1}b^{-1}g_i^{-1}(r)$ for each
$r\in\ov{j_i}$}
\end{equation}
(${\mathbf S}_2$ is the indexing set for the double sum in equation 
\eqref{s2}; terms which do not satisfy condition \eqref{m4} are zero because
they contain a degenerate simplex).
We can define a function
\[
\Upsilon:{\mathbf S}_2\to{\mathbf S}_1
\]
as follows.  Given a pair $(D,F)$ in ${\mathbf S}_2$ let $T$ denote the
quotient of $\overline{p+m-k+\sum m_i}$ by the equivalence relation $\sim$
generated by
\[
s\sim s' \mbox{ if $s$ and $s'$ are adjacent, $\phi(s)=\phi(s')$ and
$a\circ\psi(s)=a\circ\psi(s')$.}
\]
Then the diagram
\[
\xymatrix{
\{0,\ldots,p\}
&
T
\ar[l]
\ar[r]
&
\bar{m}
}
\]
where the left arrow  is induced by $\phi$ and the right by
$b\circ\psi$, is special.  Using the fact that there is a 
unique ordered bijection $T\to \overline{p+m}$ we obtain a special diagram 
\[
E:\quad\xymatrix{
\{0,\ldots,p\}
&
\overline{p+m}
\ar[l]_-{\alpha}
\ar[r]^-{\beta}
&
\bar{m}
}
\]
and this diagram satisfies \eqref{m2}.
Note
that for this diagram we have
\begin{equation}
\label{s3}
\alpha(\beta^{-1}f^{-1}(i))=\phi(\psi^{-1}b^{-1}(\ov{m_i})).
\end{equation}
Next let $T_i=\psi^{-1}b^{-1}(\overline{m_i})$.  
Using condition \eqref{m4} we see that the diagram
\[
\xymatrix{
\phi(\psi^{-1}b^{-1}\ov{m_i}) 
&
T_i
\ar[l]_-{\phi}
\ar[r]^-{b\circ\psi}
&
\ov{m_i}
}
\]
is special, 
and using equation 
\eqref{s3} we see that
\[
|T_i|= ||\beta^{-1}f^{-1}(i)|| + m_i
\]
so that we get a special diagram
\[
E_i:\quad \xymatrix{
\beta^{-1}f^{-1}(i)
&
\overline{||\beta^{-1}f^{-1}(i)||+m_i}
\ar[l]_-{\alpha_i}
\ar[r]^-{\beta_i}
&
\bar{m_i}
}
\]
and this diagram satisfies \eqref{m3}.
Finally, we define 
\[
\Upsilon(D,F)=(E,E_1,\ldots,E_k).
\]

\begin{lemma}
\label{m7}
$\Upsilon$ is a bijection.
\end{lemma}

We defer the proof to the end of this section.

Inspection of equations \eqref{revTh1} and \eqref{s2} and Definition \ref{HD} 
shows that 
Proposition \ref{gamma} will follow from:

\begin{lemma}
\label{mon1}
If $\Upsilon(D,F)=(E,E_1,\ldots,E_k)$ then
\begin{equation}
\label{mon2}
\epsilon(f,E) +
\sum_{i>i'} (m_{i'}-j_{i'})||\beta^{-1}f^{-1}(i)|| 
+\sum_i\epsilon(g_i,E_i)
\equiv \eta(D)+\epsilon(h,F)
\quad\mbox{\rm mod 2}
\end{equation}
\end{lemma}

\begin{proof}[Proof of Lemma \ref{mon1}.]
After substituting the definitions of $\epsilon(f,E)$, $\epsilon(g_i,E_i)$, 
and $\epsilon(h,F)$ (see Definition \ref{eps}) and the definition of $\eta(D)$
(see Definition \ref{HD}),
equation \eqref{mon2} becomes
\begin{multline}
\label{mon3}
\underset{f(j)>f(j')}{\underset{j<j'}\sum}
||\beta^{-1}(j)||\,||\beta^{-1}(j')||
+
\sum_{j=1}^m ||\beta^{-1}(j)||(\tau_f(j)-f(j)) 
\\
+
\sum_{i>i'} (m_{i'}-j_{i'})||\beta^{-1}f^{-1}(i)||
+
\sum_i
\underset{g_i(r)>g_i(r')}
{\underset{r<r'\text{ in $\ov{m_i}$}}\sum} 
||\beta_i^{-1}(r)|| ||\beta_i^{-1}(r')|| 
\\
+ \sum_i \sum_{r=1}^{m_i} ||\beta_i^{-1}(r)||(\tau_{g_i}(r)-g_i(r))
\quad
\equiv
\sum_{i<i'}(m_i-j_i)||f^{-1}(i')|| 
\\
+
\sum_i
\underset{g_i(r)>g_i(r')}
{\underset{r<r'\text{ in $\ov{m_i}$}}\sum}
||b^{-1}(r)|| ||b^{-1}(r')|| 
+\sum_i\sum_{r=1}^{m_i} ||b^{-1}(r)||(\tau_{g_i}(r)-g_i(r))
\\ 
+
\underset{h(s)>h(s')}{\underset{s<s'}\sum}
||\psi^{-1}(s)|| ||\psi^{-1}(s')|| +
\sum_{s=1}^{m-k+\sum m_i} ||\psi^{-1}(s)|| (\tau_h(s)-h(s))
\end{multline}

Next we substitute the equations
\[
||\beta^{-1}f^{-1}(i)||=||f^{-1}(i)||+ \sum_{f(j)=i} ||\beta^{-1}(j)||
\]
\[
||\beta^{-1}(j)||=\sum_{a(s)=j} ||\psi^{-1}(j)||
\]
and
\[
||\beta_i^{-1}(r)||=||b^{-1}(r)||+\sum_{b(s)=r}||\psi^{-1}(s)||
\]
and cancel; this reduces equation \eqref{mon3} to
\begin{multline} 
\label{mon4}
\underset{f(a(s))>f(a(s'))}{\underset{a(s)<a(s')}\sum}
||\psi^{-1}(s)|| ||\psi^{-1}(s')||
+
\sum_{s=1}^{m-k+\sum m_i} ||\psi^{-1}(s)|| (\tau_f(a(s))-f(a(s)))
\displaybreak[0]
\\
+
\sum_{f(a(s))>i'} (m_{i'}-j_{i'})||\psi^{-1}(s)||
+
\underset{g_i(b(s))>g_i(b(s'))}{\underset{b(s)<b(s')\text{ in $\ov{m_i}$}}\sum}
||\psi^{-1}(s)|| ||\psi^{-1}(s')||
\displaybreak[0]
\\
+
\underset{g_i(r)>g_i(b(s))}{\underset{r<b(s)\text{ in $\ov{m_i}$}}\sum}
||b^{-1}(r)|| ||\psi^{-1}(s)||
+
\underset{g_i(r)<g_i(b(s))}{\underset{r>b(s)\text{ in $\ov{m_i}$}}\sum}
||b^{-1}(r)|| ||\psi^{-1}(s)||
\displaybreak[0]
\\
+
\sum_{s=1}^{m-k+\sum m_i} ||\psi^{-1}(s)||(\tau_{g_i}(b(s))-g_i(b(s)))
\displaybreak[0]
\\
\equiv
\underset{h(s)>h(s')}{\underset{s<s'}\sum} ||\psi^{-1}(s)|| ||\psi^{-1}(s')||
+
\sum_{s=1}^{m-k+\sum m_i} ||\psi^{-1}(s)||(\tau_h(s)-h(s))
\end{multline} 

Now the first and fourth terms on the left side of \eqref{mon4} cancel the
second term on the right; if we 
collect the coefficients of
$||\psi^{-1}(s)||$ in all uncanceled terms, it remains to verify, for each $s$,
the following equation (where $i$ denotes $f(a(s))$)
\begin{multline}
\notag
\tau_f(a(s))-f(a(s))
+
\sum_{f(a(s))> i'} (m_{i'}-j_{i'})
+
\underset{g_i(r)>g_i(b(s))}{\underset{r<b(s)\text{ in $\ov{m_i}$}}\sum}
||b^{-1}(r)||
\\
+
\underset{g_i(r)<g_i(b(s))}{\underset{r>b(s)\text{ in $\ov{m_i}$}}\sum}
||b^{-1}(r)||
+
(\tau_{g_i}(b(s))-g_i(b(s)))
\equiv
(\tau_h(s)-h(s))
\end{multline}
and this follows from Definition \ref{tau} and the equations
\[
h(s)= g_i(s)+\sum_{i'<i} j_i 
\]
\[
\tau_h(s)=\sum_{i'<i} |b^{-1}(\ov{m_{i'}})| +|\{s'\in b^{-1}(\ov{m_i}): 
g_i(s')<g_i(s) \text{ or } 
g_i(s')=g_i(s)
\text{ and }
s'\leq s \}|
\]
\[
|b^{-1}(\ov{m_{i'}})|=m_{i'}+||f^{-1}(i')||
\]
(this equation follows from the fact that $E_{i'}$ is special)
and
\[
|\{s'\in b^{-1}(\ov{m_i}) : s'\leq s\}|=
|\{j\in f^{-1}(i):j\leq a(s)\}|
+
|\{
r\in \ov{m_i} : r\leq b(s)
\}|
-1
\]
(this equation follows from the fact that $E_i$ is special).
\end{proof}

It remains to give the proof of Lemma \ref{m7}.
We will define a map $\Phi:{\mathbf S}_1\to{\mathbf S}_2$ which is inverse to 
$\Upsilon$.
Let $E$ be a diagram of the form \eqref{E} and for each $i$ let $E_i$ be a
diagram of the form \eqref{Ei}.  Suppose that
$(E,E_1,\ldots,E_k)$ is in ${\mathbf S}_1$.  
Let $T_i=\ov{||\beta^{-1}f^{-1}(i)||+m_i}$, let 
\[
c:\coprod T_i\to \bar{m}
\]
be the map which restricts on $T_i$ to  $\alpha_i$ and let 
\[
d:\coprod T_i\to \coprod \ov{m_i}
\]
be $\coprod \beta_i$.  Give $\coprod T_i$ the unique total order for which $c$
and the inclusions of the $T_i$ are order-preserving.  Define 
\[
T=(\coprod T_i)/\sim
\]
where $\sim$ is the equivalence relation generated by
\[
\mbox{$t\sim t'$ if $t$ and $t'$ are adjacent, $c(t)=c(t')$, and $d(t)=d(t')$}
\]
Then $c$ and $d$ factor through $T$ to give a diagram
\begin{equation}
\label{m5}
\xymatrix{
T
\ar[r]^-{d'}
\ar[d]_{c'}
&
\coprod\overline{m_i}
\ar[d]_{\chi}
\\
\bar{m}
\ar[r]^{f}
&
\bar{k}
}
\end{equation}
and the diagram
\[
\xymatrix{
f^{-1}(i)
&
(d')^{-1}(\overline{m_i})
\ar[r]^-{d'}
\ar[l]_-{c'}
&
\overline{m_i}
}
\]
is special for each $i\in\bar{k}$. This implies that $|T|=m-k+\sum m_i$, so
diagram \eqref{m5} determines a diagram 
\[
D:\quad
\xymatrix{
\ov{m-k+\sum m_i}
\ar[r]^-{b}
\ar[d]_{a}
&
\coprod\overline{m_i}
\ar[d]_{\chi}
\\
\bar{m}
\ar[r]^{f}
&
\bar{k}
}
\]
which is of type $(f,m_1,\ldots,m_k)$.
Next consider the diagram
\begin{equation}
\label{m6}
\xymatrix{
\{0,\ldots,p\}
&
\coprod T_i
\ar[l]_-{\phi'}
\ar[r]^-{\psi'}
&
T
}
\end{equation}
where $\phi'$ restricts on $T_i$ to $\alpha\circ\alpha_i$ and $\psi'$ is the
quotient map.  If $q$ and $q'$ are adjacent elements of $\coprod T_i$ then
$\phi'(q)=\phi'(q')\Rightarrow \psi(q)\neq \psi(q')$ and 
$\psi'(q)=\psi'(q')\Rightarrow \phi(q)\neq \phi(q')$; since $|\coprod
T_i|=p+|T|$ this implies that \eqref{m6} is special, and using the unique
ordered bijections $\coprod T_i\to \ov{p+m-k+\sum m_i}$ and $|T|\to 
\ov{m-k+\sum m_i}$ we obtain a special diagram 
\[
F:\quad
\xymatrix{
\{0,\ldots,p\}
&
\overline{p+m-k+\sum m_i}
\ar[l]_-{\phi}
\ar[r]^-{\psi}
&
\overline{m-k+\sum m_i}
}
\]
and this diagram satisfies \eqref{m4}.
Thus the pair $(D,F)$ is in ${\mathbf S}_2$ and we may
define 
\[
\Phi(E,E_1,\ldots,E_k)=(D,F).
\]

\end{document}